\documentclass[11pt]{article}

\usepackage{amscd,amsmath, amssymb, fancyhdr, epsfig, mathbbol, dutchcal, url}
\usepackage[matrix,arrow,curve]{xy}
\usepackage[backref=page]{hyperref}
\usepackage[T1]{fontenc}

\renewcommand*{\backrefalt}[4]{%
	\ifcase #1 (Not cited.)%
	\or        (Cited on page~#2.)%
	\else      (Cited on pages~#2.)%
	\fi}

\hypersetup{
	colorlinks   = true,
	citecolor    = magenta,
	linkcolor    =blue,
	urlcolor     =magenta	
}

\numberwithin{equation}{section}

\newcommand{\version}{version 1.1,\ \ June 22, 2025}

\def\eqref#1{(\ref{#1})}

\newcommand{\g}{{\mathfrak g}}

\newcommand{\arrow}{{\:\longrightarrow\:}}
\newcommand{\Z}{{\Bbb Z}}
\def\C{{\Bbb C}}
\newcommand{\R}{{\Bbb R}}

\renewcommand{\H}{{\Bbb H}}

\def\1{\sqrt{-1}\:}
\newcommand{\restrict}[1]{{\left|_{{\phantom{|}\!\!}_{#1}}\right.}}
\newcommand{\cntrct}                % contraction with a vector field
{\hspace{2pt}\raisebox{1pt}{\text{$\lrcorner$}}\hspace{2pt}}
\newcommand{\calo}{\mathcal O}

\renewcommand{\phi}{\varphi}
\renewcommand{\epsilon}{\varepsilon}

\newcommand{\End}{\operatorname{End}}

\newcommand{\Tot}{\operatorname{Tot}}
\newcommand{\Id}{\operatorname{Id}}

\newcommand{\Vol}{{\operatorname{Vol}}}

\newcommand{\Sec}{\operatorname{Sec}}
\newcommand{\Hom}{\operatorname{Hom}}
\newcommand{\Sym}{\operatorname{Sym}}
\newcommand{\Hol}{\operatorname{Hol}}

\newcommand{\Aff}{\operatorname{Aff}}

\newcommand{\Hor}{\operatorname{Hor}}

\newcommand{\Tw}{\operatorname{Tw}}

\newcommand{\Alb}{\operatorname{Alb}}

\newcounter{Mycounter}[section]
\newcounter{lemma}[section]
\renewcommand{\thelemma}{{Lemma \thesection.\arabic{lemma}}}
\newcommand{\lemma}{%
    \setcounter{lemma}{\value{Mycounter}}
    \refstepcounter{lemma}
    \stepcounter{Mycounter}
    {\noindent \bf \thelemma:\ }}

\newcounter{claim}[section]
\renewcommand{\theclaim}{{Claim \thesection.\arabic{claim}}}
\newcommand{\claim}{%
    \setcounter{claim}{\value{Mycounter}}
    \refstepcounter{claim}
    \stepcounter{Mycounter}
    {\noindent \bf \theclaim:\ }}

\newcounter{sublemma}[section]
\newcounter{corollary}[section]
\renewcommand{\thecorollary}{{Corollary \thesection.\arabic{corollary}}}
\newcommand{\corollary}{%
    \setcounter{corollary}{\value{Mycounter}}
    \refstepcounter{corollary}
    \stepcounter{Mycounter}
    {\noindent \bf \thecorollary:\ }}

\newcounter{theorem}[section]
\renewcommand{\thetheorem}{{Theorem \thesection.\arabic{theorem}}}
\newcommand{\theorem}{%
    \setcounter{theorem}{\value{Mycounter}}
    \refstepcounter{theorem}
    \stepcounter{Mycounter}
    {\noindent \bf \thetheorem:\ }}

\newcounter{conjecture}[section]
\renewcommand{\theconjecture}{{Conjecture \thesection.\arabic{conjecture}}}
\newcommand{\conjecture}{%
    \setcounter{conjecture}{\value{Mycounter}}
    \refstepcounter{conjecture}
    \stepcounter{Mycounter}
    {\noindent \bf \theconjecture:\ }}

\newcounter{proposition}[section]
\renewcommand{\theproposition}
      {{Proposition \thesection.\arabic{proposition}}}
\newcommand{\proposition}{%
    \setcounter{proposition}{\value{Mycounter}}
    \refstepcounter{proposition}
    \stepcounter{Mycounter}
    {\noindent \bf \theproposition:\ }}

\newcounter{definition}[section]
\renewcommand{\thedefinition}
      {{Definition~\thesection.\arabic{definition}}}
\newcommand{\definition}{%
    \setcounter{definition}{\value{Mycounter}}
    \refstepcounter{definition}
    \stepcounter{Mycounter}
    {\noindent \bf \thedefinition:\ }}

\newcounter{example}[section]
\renewcommand{\theexample}{{Example \thesection.\arabic{example}}}
\newcommand{\example}{%
    \setcounter{example}{\value{Mycounter}}
    \refstepcounter{example}
    \stepcounter{Mycounter}
    {\noindent \bf \theexample:\ }}

\newcounter{remark}[section]
\renewcommand{\theremark}{{Remark \thesection.\arabic{remark}}}
\newcommand{\remark}{%
    \setcounter{remark}{\value{Mycounter}}
    \refstepcounter{remark}
    \stepcounter{Mycounter}
    {\noindent \bf \theremark:\ }}

\newcounter{problem}[section]
\newcounter{question}[section]
\newcommand{\proof}{\noindent{\bf Proof:\ }}

\newcommand{\pstep}{\noindent{\bf Proof. Step 1:\ }}

\makeatletter

\renewcommand{\leftmark}{{\scriptsize Exotic hypercomplex
    structures on a torus}}

\@addtoreset{equation}{section} \@addtoreset{footnote}{section}

\def\x@arrow{\DOTSB\Relbar}
\def\xlongrightarrowfill@{\arrowfill@\relbar\relbar\longrightarrow}
\newcommand{\xlongrightarrow}[2][]{%
        \ext@arrow 0099\xlongrightarrowfill@{#1}{#2}}
\makeatother

\def\blacksquare{\hbox{\vrule width 5pt height 5pt depth 0pt}}
\def\endproof{\blacksquare}

\begin{document}

%%%%%%%%%%%%%%%%%%%%%%%%%%%%%%%%%%%%%%%%%%%%%%%%%%%%%%%%%%%%
\begin{center}
{\LARGE\bf
Exotic hypercomplex structures on a torus do not exist\\[4mm]
}
%%%%%%%%%%%%%%%%%%%%%%%%%%%%%%%%%%%%%%%%%%%%%%%%%%%%%%%%%%%%%

Alberto Pipitone Federico, 
Misha Verbitsky\footnote{M. V. acknowledges support 
of CNPq - Process 310952/2021-2,
and FAPERJ SEI-260003/000410/2023. \\

}

\end{center}

{\small \hspace{0.10\linewidth}
\begin{minipage}[t]{0.85\linewidth}
{\bf Abstract.} 
A hypercomplex manifold is a manifold with three complex
structures satisfying quaternionic relations. 
Such a manifold admits a unique torsion-free connection
preserving the quaternionic action, called {\bf the Obata connection.}
A compact K\"ahler manifold admitting a hypercomplex
structure always admits a hyperk\"ahler structure as well;
however, it is not obvious whether the original
hypercomplex structure is hyperk\"ahler. 
A non-hyperk\"ahler hypercomplex structure on a
K\"ahler manifold is called {\bf exotic}. 
We show that the Obata connection for an 
exotic hypercomplex structure on a torus is flat
and classify complete flat affine structures
on real tori. We use this classification
to prove that exotic hypercomplex structures
do not exist.
\end{minipage}
}

\tableofcontents

%%%%%%%%%%%%%%%%%%%%%%%%%%%%%%%%%%%%%%%%%%%%%%%%%%%%%%%%%%%%%%%%%%%%%%%%%%%%

\section{Introduction}

%%%%%%%%%%%%%%%%%%%%%%%%%%%%%%%%%%%%%%%%%%%%%%%%%%%%%%%%%%%%%%%%%%%%%%%%%%%%

Throughout this paper, a ``complex manifold''
is an almost complex manifold $(M,I)$, with the tensor 
$I \in \End(TM)$, $I^2=-\Id$ being integrable.
In this situation, the operator $I$ is called 
``the complex structure''.

A hypercomplex manifold is a manifold $M$
equipped with three complex structures $I, J, K$
satisfying the quaternionic relations
$IJ=-JI=K$. This defines, in particular, a
quaternionic action on the tangent bundle,
hence the real dimension of $M$ is divisible by four.

As shown by Obata (\cite{_Obata_}),
every hypercomplex manifold admits a unique torsion-free
connection preserving the quaternionic action, known as 
{\bf Obata connection}; its holonomy group is
contained in $GL(n, {\Bbb H})$. Conversely,
if a manifold admits a torsion-free connection
with holonomy in $GL(n, {\Bbb H})$, it 
admits three complex structures $I, J, K$
preserved by the connection and
satisfying the quaternionic relations.
A complex structure preserved by
a torsion-free connection is always
integrable (\cite[Example I.2.4]{_Kobayashi:transformations_},
see also \cite{_Frolicher:complex_}), hence
hypercomplex manifolds can be
defined as manifolds equipped with 
a torsion-free connection with 
holonomy in $GL(n, {\Bbb H})$.

Hyperk\"ahler manifolds
are  manifolds equipped with a hypercomplex structure $I,J, K$
and a Riemannian metric which is K\"ahler with respect
to $I, J , K$. However, this notion also has an equivalent
definition in terms of holonomy: a Riemannian manifold
is hyperk\"ahler if the holonomy of its Levi-Civita
connection belongs to $Sp(n)$ (\cite{_Besse:Einst_Manifo_}).
Clearly, this is equivalent to holonomy of Obata connection
preserving a metric, or having a compact closure.

By definition, the Levi-Civita connection on a Riemannian
manifold is the unique torsion-free connection preserving
the metric, thus if the Obata connection preserves a metric, 
it is automatically the Levi-Civita connection of this
metric. In this case its holonomy belongs
to $O(4n) \cap GL(n, {\Bbb H})=Sp(n)$.
In other words, a hypercomplex manifold is hyperk\"ahler
if and only if its Obata connection preserves a metric.

In this paper we study the compact hypercomplex manifolds
$(M,I,J,K)$ such that $(M,I)$ admits a K\"ahler structure.
As shown in \cite{_Verbitsky:Kahler_HKT_}, in this case
$(M,I)$ admits a hyperk\"ahler structure, i.e.
there is a hypercomplex structure $(I,J',K')$ on $M$
such that its Obata holonomy belongs to $Sp(n)$.

However, this result does not imply that $(I, J, K)$ 
is hyperk\"ahler. We call a hypercomplex structure $(M,I,J,K)$
{\bf exotic} if the holonomy of its Obata connection
does not preserve any metric, but $(M,I)$ admits a K\"ahler structure.

In \cite{_Verbitsky:Kahler_HKT_} it was conjectured 
that exotic hypercomplex structures does not exist.
If $\dim_\R M=4$, this statement is contained
\cite{_Boyer:Hyperhermitian_}.\footnote{
For the reader's convenience, we 
give a proof of non-existence of exotic
hypercomplex structures on K3 in 
Subsection \ref{_hc_K3_Subsection_}.}

In the present paper, we prove it for 
complex tori of arbitrary dimension. 

Note that not all complex structures on
a manifold diffeomorphic to a compact torus are of K\"ahler type.
The first examples of complex structures of non-K\"ahler type
on a torus were obtained by A. Blanchard (\cite{_Blanchard:Recherche_})
and independently by A. Sommese (\cite{_Sommese:quaternionic_}).
An introduction to the history and current research on
the subject is found in this paper of F. Catanese,
\cite{_Catanese:Deformation_in_the_large}; see also
\ref{_twi_double_cover_Example_} where one of such
examples is constructed. 

In the present paper, we define a {complex torus}
as a compact quotient $\C^{n}/\Z^{2n}$, where $\Z^{2n}$
acts by translations; this quotient is clearly of K\"ahler type.

A {\bf flat affine manifold} is a manifold equipped
with a torsion-free flat connection; it is {\bf complete}
if it obtained as a quotient of $\R^n$ by an affine action
of a discrete group.
For an introduction to flat affine geometry,
see Subsection \ref{_flat_affine_Subsection_}. 

The first step towards our main result 
is the following theorem.

\hfill

\theorem
Let $(M,I,J,K)$ be a hypercomplex manifold.
Assume that $(M,I)$ is biholomorphic to a complex
torus. Then the Obata connection of $(M,I,J,K)$
is flat and complete.

\proof \ref{_exotic_Obata_flat_Theorem_}. \endproof

\hfill

In Section \ref{_flat_affine_on_torus_Section_} 
we classify complete flat affine structures on 
manifolds diffeomorphic to a torus. Surprisingly,
this classification is non-trivial. Fortunately
for our research, the key results were obtained
in 1981 by Fried, Goldman and Hirsch
(\ref{_FGH_Theorem_}). Using this classification, 
we were able to produce examples of non-trivial flat
affine structures on complex manifolds.
However, the same construction, when applied
to hypercomplex manifolds, leads to a
contradiction, proving our main result:

\hfill

\theorem
Exotic hypercomplex structures on complex tori
do not exist.

\proof
\ref{_no_hc_tori_exotic_Theorem_}.
\endproof

%%%%%%%%%%%%%%%%%%%%%%%%%%%%%%%%%%%%%%%%%%%%%%%%%%%%%%%%%%%%%%%%%%%%%%%%

\section{Hypercomplex manifolds}

%%%%%%%%%%%%%%%%%%%%%%%%%%%%%%%%%%%%%%%%%%%%%%%%%%%%%%%%%%%%%%%%%%%%%%%%

%%%%%%%%%%%%%%%%%%%%%%%%%%%%%%%%%%%%%%%%%%%%%%%%%%%%%%%%%%%%%%%%%%%%%%%%
\subsection{Hypercomplex structures and Obata connection}
%%%%%%%%%%%%%%%%%%%%%%%%%%%%%%%%%%%%%%%%%%%%%%%%%%%%%%%%%%%%%%%%%%%%%%%%

\definition Let $M$ be a smooth
manifold equipped with endomorphisms
$I, J, K:\; TM\arrow TM$, satisfying the quaternionic relation
$I^2=J^2=K^2=IJK=-\Id.$  Suppose that $I$, $J$, $K$ are
integrable. Then $(M, I, J, K)$ is called {\bf a hypercomplex manifold}.

\hfill

\theorem  Let $(M,I,J,K)$ be 
a hypercomplex manifold. Then $M$ admits a unique torsion-free
affine connection preserving $I, J, K$. \\
\proof \cite{_Obata_}. 
\endproof

\hfill

\definition
This connection is called {\bf the Obata connection}.

\hfill

\remark
Holonomy of Obata connection lies in $GL(n, {\Bbb H})$.
Conversely, a manifold equipped with an affine,
torsion-free connection with holonomy in $GL(n, {\Bbb H})$ 
is hypercomplex. Indeed, suppose that $I,J,K$ are operators defining quaternionic structure
on $TM$ and  $\nabla$ is an affine torsion-free connection 
preserving $I$, $J$, $K$.  Then $I$, $J$, $K$ are integrable
almost complex structures and $(M,I,J,K)$ is hypercomplex. 

\hfill

\remark
This can be used as a definition of a hypercomplex 
structure: {\bf  a hypercomplex manifold} $(M, \nabla, I, J, K)$ is a 
manifold equipped with a torsion-free connection such that
its holonomy preserves a quaternionic structure on a tangent bundle.

\hfill

\remark Compatibility of $\nabla$ with $I,J,K$ is equivalent to having $\nabla^{0,1}_I= \overline{\partial }_I$, $\nabla^{0,1}_J= \overline{\partial }_J $, $\nabla^{0,1}_K= \overline{\partial }_K$. The curvature form of Obata's connection is of type $(1,1)$ with respect to $I,J,K$. This follows directly from Bianchi identity, the fact that $\nabla$ preserves $I,J,K$ and that $J$ maps $T^{1,0}_I M$ to $T^{0,1}_I M$.

%%%%%%%%%%%%%%%%%%%%%%%%%%%%%%%%%%%%%%%%%%%%%%%%%%%%%%%%%%%%%%%%%%%%%%%%
\subsection{Hyperk\"ahler structures}
%%%%%%%%%%%%%%%%%%%%%%%%%%%%%%%%%%%%%%%%%%%%%%%%%%%%%%%%%%%%%%%%%%%%%%%%

The following definition is equivalent to the standard one;
we leave the equivalence as an exercise to the reader.

\hfill

\definition
A hypercomplex manifold $(M, \nabla, I, J, K)$ is called
{\bf hyperk\"ahler} if the holonomy $\Hol(\nabla)$ of
$\nabla$ preserves a 
quaternionic invariant Riemannian metric $g$. Such a metric
is called {\bf  hyperk\"ahler}. A {\bf hyperk\"ahler structure}
is $(M, \nabla, I, J, K, g)$. In this situation,
 $\nabla$ is the Levi-Civita connection.

\hfill

\theorem
Let $(M,I,J,K)$ be a compact hypercomplex manifold.
Assume that $(M,I)$ admits a K\"ahler structure.
Then $(M,I)$ admits a hyperk\"ahler structure
$(I,J',K')$.

\proof \cite{_Verbitsky:Kahler_HKT_}. \endproof

\hfill

\definition
Let $(M,I,J,K)$ be a compact hypercomplex manifold.
Assume that $(M,I)$ admits a K\"ahler structure.
The hypercomplex structure $(I,J,K)$ is called
{\bf exotic} if it is not compatible
with a hyperk\"ahler metric, that is, if the
holonomy of its Obata connection is non-compact.

%%%%%%%%%%%%%%%%%%%%%%%%%%%%%%%%%%%%%%%%%%%%%%%%%%%%%%%%%%%%%%%%%%%%%%%%
\subsection{Exotic hypercomplex structures on K3}
\label{_hc_K3_Subsection_}
%%%%%%%%%%%%%%%%%%%%%%%%%%%%%%%%%%%%%%%%%%%%%%%%%%%%%%%%%%%%%%%%%%%%%%%%

The following, as well as the same
result for two dimensional complex tori, can be found in 
\cite{_Boyer:Hyperhermitian_}
We give a short proof of this result for the benefit of the reader. 

\hfill

\theorem
 Exotic hypercomplex structures on a K3 surface do not exist. 

\hfill

\pstep
Let $K_{M_I}$, $K_{M_J}$, $K_{M_K}$ denote the canonical
bundle of a hypercomplex manifold $M$ with the complex structures $I, J, K$.
Since the action of $SU(2)=U({\Bbb H}, 1)$ on $TM$ interchanges
the complex structures $I, J,K$, it also interchanges
$K_{M_I}$, $K_{M_J}$, $K_{M_K}$ considered as sub-bundles
in $\Lambda^{2n}(M)$. This implies that $K_{M_I}$, $K_{M_J}$, $K_{M_K}$
are isomorphic as complex vector bundles with the connection
induced by the Obata connection.

Let $(M,I,J,K)$ be a hypercomplex structure on 
a K3 surface and $\Theta$ the curvature of the Obata connection on 
its canonical bundle $K_{M_I}=K_{M_J}=K_{M_K}$. Since $\Theta$ is a $(1,1)$-form with respect to $I,J,K$, it is $SU(2)$-invariant with respect to the
$SU(2)$-action on $\Lambda^*(M)$ generated by quaternions. 
In particular, for any Hermitian metric $g$, one has $* \Theta=-\Theta$, thus
\[
\Theta\wedge \Theta = -\|\Theta\|^2_g \Vol_g
\]
By Bianchi identity,  $\Theta$ is closed, and cohomologous
to $c_1(M,I)$. Since $(M,I)$ admits a non-degenerate
$(2,0)$-form,  we have $c_1(M,I)=0$. 
Integrating the above formula, we find $\Theta=0$.

Given that $\pi_1(K3)=0$,  we obtain that $K_{M_I}$ (resp. $K_{M_J}$, $K_{M_K}$)
is holomorphically trivialized by an Obata-parallel section $s_I$ (resp. $s_J$, $s_K$).

\hfill

{\bf Step 2:}
The holomorphic sections $s_I, s_J, s_K$
generate a rank 3 sub-bundle $B$ in $\Lambda^2(M)$,
which has flat connection because $s_I, s_J, s_K$ are parallel.
Recall that the natural map $SL(4,\R)\to SO(3,3)$ taking
endomorphisms of $\R^2$ to endomorphisms of $\Lambda^2(\R^4)$
is a 2 to 1 covering. Fixing a maximal compact subgroup
$SO(3)\times SO(3)$ in $SO(3,3)$ is the same as fixing
a metric on $\R^4$; this reduces the structure group
$SL(4, \R)$ to $SO(4)\subset  SL(4, \R)$, which is also
a 2:1 covering of $SO(3)\times SO(3)$. Finally, reducing
the structure group of $\Lambda^2(\R^4)$ to $SO(3)\subset SO(3)\times SO(3)$
means that the structure group $SL(4, \R)$ is reduced to $Sp(1)=SU(2)$,
which is a 2:1 cover of $SO(3)$. In our case, this reduction is
given by a trivialization of $B$. Therefore, the holonomy
of $\nabla$ is contained in $Sp(1)$, hence the hypercomplex
structure on $M$ is in fact hyperk\"ahler.
\endproof

%%%%%%%%%%%%%%%%%%%%%%%%%%%%%%%%%%%%%%%%%%%%%%%%%%%%%%%%%%%%%%%%%%%%%%%%

\section{Twistor spaces for hypercomplex manifolds}

%%%%%%%%%%%%%%%%%%%%%%%%%%%%%%%%%%%%%%%%%%%%%%%%%%%%%%%%%%%%%%%%%%%%%%%%

%%%%%%%%%%%%%%%%%%%%%%%%%%%%%%%%%%%%%%%%%%%%%%%%%%%%%%%%%%%%
\subsection{Twistor spaces: definition and basic notions}
%%%%%%%%%%%%%%%%%%%%%%%%%%%%%%%%%%%%%%%%%%%%%%%%%%%%%%%%%%%%

\definition 
{\bf Induced complex structures} 
on a hypercomplex manifold are 
complex structures of the form $L= aI + bJ +c K$, with $a^2+b^2+c^2=1$.

\hfill

\remark
The induced complex structures are usually non-algebraic. Indeed,
if $M$ is compact and hyperk\"ahler, for generic $a, b, c$, $(M,L)$ has no divisors
(\cite{_Fujiki_}). 

\hfill

A twistor space $\Tw(M)$ of a hypercomplex manifold
is a complex manifold obtained by gluing 
these complex structures into
a holomorphic family over $\C P^1$. More formally:

\hfill

\definition
 Let $\Tw(M)$ be the smooth manifold
$M \times S^2$. Consider the complex structure $I_m:T_mM \to T_mM$ 
on $M$ induced by $J \in S^2 \subset {\Bbb H}$. Let $I_J$
denote the complex structure on $S^2 = \C P^1$.
The operator 
\[
I_{\Tw} = I_m \oplus I_J:T_x\Tw(M) \to T_x\Tw(M)
\]
satisfies $I_{\Tw} ^2 = -\Id$, defining
an almost complex structure on $\Tw(M)$ which is known to be integrable
(\cite{_Kaledin_}).

\hfill

\example If $M={\Bbb H}^n$, then 
$\Tw(M)= \Tot (\calo(1)^{\oplus n}) \cong \C P^{2n+1} \backslash \C P^{2n-1}$
(see e. g. \cite{_JV:Instantons_}).

\hfill

\remark\label{_O(1)_normal_Remark_}
A curve $\C P^1\times \{m\}\subset \Tw(M)$ is by
construction holomorphic; such a curve is called
{\bf a horizontal curve}. Its normal bundle
is $\calo(1)^{\oplus n}$ (\cite{_HKLR_}).

\hfill

\remark 
For $M$ compact,  $\Tw(M)$ never admits a K\"ahler, or Moishezon structure
(\cite{_Pipitone:twistor_,_Gorginyan:Moishezon_}).

%%%%%%%%%%%%%%%%%%%%%%%%%%%%%%%%%%%%%%%%%%%%%%%%%%%%%%%%%%%%%
%\subsection{Rational curves on $\Tw(M)$}
%%%%%%%%%%%%%%%%%%%%%%%%%%%%%%%%%%%%%%%%%%%%%%%%%%%%%%%%%%%%%
%
%\definition
%{\bf An ample rational curve} on a complex manifold $M$ is
%a smooth curve $S \cong \C P^1\subset M$ such that 
%$NS=\bigoplus_{k=1}^{n-1}\calo(i_k)$, with $i_k >0$.
%It is called {\bf  a quasiline} if all
%$i_k=1$.
%
%\hfill
%
%\theorem 
%Let $M$ be a compact complex manifold containing a 
%an ample rational line. Then any $N$ points $z_1, ..., z_N$ can
%be connected by an ample rational curve.
%
%\proof
%\cite{_V:curves_twistor_}.
%\endproof
%
%\hfill
%
%\claim
%Let $M$ be a hyperk\"ahler
%manifold, $\Tw(M)\stackrel \sigma\arrow M$ its twistor
%space, $m\in M$ a point, and $S_m=\C P^1 \times \{m\}$
%the corresponding rational curve in $\Tw(M)$.  Then $S_m$
%  is a quasiline.
%
%\hfill
%
%\proof Since the claim is essentially
%infinitesimal, it suffices to check it when $M$ is flat.
%Then $\Tw(M)= \Tot (\calo(1)^{\oplus 2p})
%\cong \C P^{2p+1} \backslash \C P^{2p-1}$, and
%$S_m$ is a section of $\calo(1)^{\oplus 2p}$.
%\endproof

%%%%%%%%%%%%%%%%%%%%%%%%%%%%%%%%%%%%%%%%%%%%%%%%%%%%%%%%%%%%
\subsection{The twistor data}
%%%%%%%%%%%%%%%%%%%%%%%%%%%%%%%%%%%%%%%%%%%%%%%%%%%%%%%%%%%%

Hypercomplex structures can be characterized (and,
in fact, defined) in terms of twistor spaces.
The main advance in this direction was made
by Hitchin, Karlhede, Lindstr\"om, Ro\v cek in \cite{_HKLR_};
see also \cite{_Simpson:hyperka-defi_,_Verbitsky:hypercomple_}.
We use the notion of ``the twistor data''
of a hypercomplex manifold to signify the twistor space
equipped with additional structures needed to reconstruct
the hypercomplex manifold.

\hfill

Let $\check \tau$ denote the central symmetry
on $\C P^1$; if we identify $\C P^1$ with imaginary unit
quaternions, we have
$\check\tau(L)=-L$. It is an anticomplex involution without fixed points.

\hfill

\definition\label{_twistor_data_Definition_}
The {\bf twistor data}
is a complex manifold $\Tw$, $\dim \Tw=2n+1$, equipped with the following
structures. 
\begin{description}
\item[(i)]
A holomorphic submersion $\pi:\; \Tw\arrow \C P^1$
and an anticomplex involution $\tau:\; \Tw\arrow \Tw$
which makes this diagram commutative:
\[\begin{CD} \Tw @>\tau >> \Tw\\
@V\pi VV @VV\pi V\\
\C P^1 @>{\check \tau}>> \C P^1 
\end{CD}
\]
\item[(ii)]
 A connected component $\Hor$
in the set $\Sec^\tau \subset \Sec$ of $\tau$-invariant sections of $\pi$
such that for each $S\in \Hor$, the normal bundle to $S$
is $\calo(1)^{2n}$ and for each point $x\in \Tw$ there exists a unique $S\in \Hor$ passing through $x$.
\end{description}

\hfill

\remark
With any twistor space $\Tw(M)$ of a hypercomplex manifold,
one associates the twistor data in a natural way:
$\tau(I, m)= (-I, m)$, and $\Hor(M)$ the space
of all sections $S_m$, taking $I\in \C P^1$ to
$(I,m)\in \Tw(M)$, where $m\in M$ is a fixed point.

\hfill

\theorem \label{_twistor_data_determine_hc_Theorem_}
Let $M$ be a hypercomplex manifold. Then 
the twistor data on $\Tw(M)$ can be used to recover the
hypercomplex structure on $M$, which is identified with $\Hor$
by \ref{_twistor_data_Definition_} (ii).
Moreover, for any twistor data $(\Tw, \tau, \Hor)$, 
there exists a hypercomplex structure $(I,J,K)$ on $\Hor$ 
such that these twistor data are associated with
$(I,J,K)$.

\proof
\cite{_HKLR_}. 
\endproof

\hfill

\remark
This theorem has a version appropriate for hyperk\"ahler
metrics (also in \cite{_HKLR_}).

%%%%%%%%%%%%%%%%%%%%%%%%%%%%%%%%%%%%%%%%%%%%%%%%%%%%%%%%%%%%%%%%%%%%%%%%

\section{Flat affine structures and Obata connection}

%%%%%%%%%%%%%%%%%%%%%%%%%%%%%%%%%%%%%%%%%%%%%%%%%%%%%%%%%%%%%%%%%%%%%%%%

%%%%%%%%%%%%%%%%%%%%%%%%%%%%%%%%%%%%%%%%%%%%%%%%%%%%%%%%%%%%
\subsection{Flat affine manifolds}
\label{_flat_affine_Subsection_}
%%%%%%%%%%%%%%%%%%%%%%%%%%%%%%%%%%%%%%%%%%%%%%%%%%%%%%%%%%%%

To approach the exotic hypercomplex structures on a torus, we
would introduce flat affine manifolds. The following 
notions are standard (see
e. g. \cite{_Abels:survey_}, \cite{_FGH:Affine_} or
 \cite[Chapter 27]{_OV:Principles_}).

\hfill

\definition
A {\bf flat affine structure} on a manifold $M$
is a flat torsion-free connection on the tangent bundle.

\hfill

It can be easily seen that this is equivalent to the existence of an affine atlas. \\

\definition\label{affine map}
A smooth map between smooth manifolds $f:M\to N$ is
{\bf compatible with the affine structures}, or {\bf affine} if the differential
$Df$ takes parallel vector fields to parallel vector fields.

\hfill

\remark Equivalently, $f$ is affine if local parallel
forms on $N$ are pulled back to parallel forms on $M$. This
is happens precisely if there exist affine charts such
that $f$ in local coordinates is an affine map. Any
diffeomorphism compatible with the affine structures is an
affine isomorphism.

\hfill

\definition
Let $M$ be a simply connected affine manifold,
and $\theta_1, ..., \theta_n \in \Lambda^1 M$
a basis of parallel 1-forms. Since a parallel 1-form
is closed and $H^1(M,\R)=0$, the forms $\theta_i$ are
exact. Then $\theta_i = d x_i$. 
The map $\delta:\; M\to \R^n$ taking $m$ to $(x_1(m), ..., x_n(m))$
is called {\bf the development map.} We consider
$\R^n$ as a flat affine manifold, with the standard
flat affine structure. The development map
is compatible with the affine structure:

\hfill

\claim \label{_compatible_aff_Claim_}
The development map $\delta:\; M\to \R^n$
 is compatible with the flat affine connections.

\hfill

\proof It takes the coordinate 1-forms $dx_1, ...,
dx_n\in \Lambda^1(M)$
to $\theta_1, ..., \theta_n \in \Lambda^1 M$.
However, these 1-forms are parallel.
\endproof

\hfill

\definition
A flat $n$-dimensional affine manifold $(M, \nabla)$ is said
{\bf complete} if its universal cover is isomorphic to
$\R^n$ as a flat affine manifold. In this case,
$\pi_1(M)\to \text{Diff}(\R^n)$ factors through
$\Aff(\R^n)$.\footnote{The image of $\pi_1(M)\to
  \text{Diff}(\R^n)$ 
is identified with the affine holonomy,
see \ref{_affine_holo_Definition_} below.}
Equivalently, $M$ is complete if and only if $\nabla$ is geodesically complete or if the developing map is a diffeomorphism.

\hfill

\remark
Marcus conjecture (\ref{_Marcus_conjecture_}) claims that a compact flat affine manifold
is complete if and only if the flat connection preserves a volume form.

%%%%%%%%%%%%%%%%%%%%%%%%%%%%%%%%%%%%%%%%%%%%%%%%%%%%%%%%%%%%
\subsection{Linear and affine holonomy}
%%%%%%%%%%%%%%%%%%%%%%%%%%%%%%%%%%%%%%%%%%%%%%%%%%%%%%%%%%%%

\definition
The {\bf linear holonomy} (or {\bf holonomy})
of a flat affine connection $\nabla$ is its monodromy in $TM$;
by definition, the holonomy group belongs to $GL(T_x M)$,
where $x\in M$ is a base point. 

%The holonomy of an affine manifold is the holonomy of the
%connection induced by the affine structure.

\hfill

\definition
Let $\Aff(\R^n)$ denote the group of affine transforms of
$\R^n$. Clearly, $\Aff(\R^n)$ is a semidirect product,
$\Aff(\R^n)= GL(n, \R)\rtimes \R^n$.
The natural map $\Aff(\R^n)\arrow GL(n, \R)$
is called {\bf the linearization.}

\hfill

\definition\label{_affine_holo_Definition_}
Let $M$ be a flat affine $n$-manifold,
$\tilde M$ its universal cover\footnote{Any finite cover of an affine manifold admits a canonical affine structure.}  and
$\delta:\; \tilde M\to \R^n$ the development map. For any $\gamma \in \pi_1(M)$, let $\tilde{\gamma}$ be the induced homeomorphism $\tilde{M} \to \tilde{M}$.
Then define ${a}:\; \pi_1(M) \arrow \Aff(\R^n)$ to be the map making the following diagram commute
\[\begin{CD}
\tilde M @>\delta >> \R^n\\
@V\tilde{\gamma} VV  @VV{ a(\gamma)}V\\
\tilde M @>\delta >> \R^n.
\end{CD}\]
The map ${\cal a}:\; \pi_1(M) \arrow \Aff(\R^n)$ 
is called {\bf affine holonomy map.} 
In terms of transition functions, the map $a$ can be described explicitly: fix $x_0\in M$ and let $\gamma$ be a smooth loop based in $x_0$ and contained in the affine charts $(U_i,
\phi_i)_{i=1}^k$, where $x_0\in U_1 \cap U_k$ and every $U_i\cap U_{i+1}$ is non-empty and connected. Then $a(\gamma)= \phi_1\phi_k\phi_k^{-1}
\phi_{k-1} \dots \phi_2 \phi_1^{-1}$.

\hfill

For affine manifolds, the two notions are related by the following

\hfill

\lemma
The holonomy of an affine manifold is 
the linearization if its affine holonomy.

\hfill

\proof
It is clear from the above description of $a$ in terms of transition functions.
\endproof

%%%%%%%%%%%%%%%%%%%%%%%%%%%%%%%%%%%%%%%%%%%%%%%%%%%%%%%%%%%%
\subsection{Complex tori: definition and basic notions}
%%%%%%%%%%%%%%%%%%%%%%%%%%%%%%%%%%%%%%%%%%%%%%%%%%%%%%%%%%%%

\definition
A {\bf complex torus} is a complex manifold $M$ biholomoprhic to $\C^n/ \Lambda$, where $\Lambda$ is a $2n$ dimensional integral lattice. In particular, $M$ is K\"ahler and the Albanese map 
\[
\Alb: M \arrow\frac{H^{0}(M,\Omega^1)^*}{H_1(M, \Z)}
\]
is an isomorphism. 

\hfill

\remark
Any K\"ahler manifold $X$
diffeomorphic to a torus has this nature: in
(\cite{_Catanese:Deformation_in_the_large}) it is indeed
proven that, in the class of K\"{a}hler manifolds, being a
complex torus is equivalent to the purely topological
condition $\Lambda H^1(X,\Z) \backsimeq H^{\bullet}(X,\Z)$
(the isomorphism is in the category of rings).

\hfill

Non-K\"ahler complex structures on a manifold
diffeomorphic to a torus were constructed in
\cite{_Blanchard:Recherche_,_Sommese:quaternionic_}.
For the convenience of the reader, we include
one of such constructions here (see also
\cite{_Catanese:Deformation_in_the_large}).

\hfill

\example\label{_twi_double_cover_Example_}
Let $M$ be a hyperk\"ahler torus, and
$\Tw(M)$ its twistor space, equipped with
a holomorphic projection $\Tw(M) \to \C P^1$.
Let  $E$ be an elliptic
curve equipped with the standard 2-sheeted
ramified cover map $\sigma: E \to \C P^1$.
Since $\Tw(M)$ is diffeomorphic to a 
product $M \times \C P^1$, the fibered
product $E\times_{\C P^1} \Tw(M)$ is
diffeomorphic to a torus.
To see that $E\times_{\C P^1} \Tw(M)$ 
is not of K\"ahler type, we use the formula \cite[(8.12)]{_NHYM_}:
 $dd^c \omega=\omega\wedge \pi^*\omega_{\C P^1}$,
where $\omega$ is the standard Hermitian form on 
the twistor space, and $\omega_{\C P^1}$ the Fubini-Study
form on $\C P^1$. This form is clearly positive
and exact, but a K\"ahler manifold cannot
admit positive exact forms (otherwise this
form, integrated with an appropriate power
of the K\"ahler form, would be positive, which
is impossible because it is exact). 
Since $E\times_{\C P^1} \Tw(M)$ 
is a ramified double cover of the twistor space,
this manifold also admits a positive exact form,
hence it cannot be K\"ahler.

\hfill

Further on in this paper, we will use the following
important theorem.

\hfill

\theorem  
Let ${\cal X}$ be a connected family of complex
structures on a manifold $M$ diffeomorphic to a torus.
Assume that for some $I\in {\cal X}$, the manifold $(M,I)$
is a complex torus.  Then $(M,I')$ is a torus for all
$I'\in {\cal X}$.

\proof \cite[Theorem 4.1]{_Catanese:Deformation_types_}. \endproof

\hfill

%\definition
%A {\bf holomorphic connection} on a complex manifold
%is a holomorphic differential
% operator $\nabla:\; {\cal T} \arrow {\cal T} \otimes_{\calo_M}\Omega^1 M$ 
%acting on the sheaf ${\cal T}$of holomorphic vector fields which satisfies the Leibniz rule
%$\nabla(f X)= df\otimes X + f\nabla(X)$
%
%\hfill

\remark\label{_flat_on_torus_Remark_}
Let $\theta_1, ... \theta_n$ be holomorphic
differentials on a complex torus $M$. Their antiderivatives
define a flat affine chart on $M$; the corresponding
flat affine structure on $M$ is canonically defined.
This also defines a flat affine 
connection on $M$ compatible with the holomorphic structure.

%\hfill
%
%\remark Also, each complex torus $M$
%is a torsor over the corresponding group manifold,
%identified with the connected component $\Aut_0(M)$ of $\Aut(M)$,
%and its action on $M$ is canonically defined. Since
%$\Aut_0(M)$ is (non-canonically) identified with $M$,
%this action is called {\bf the action of the
%torus on itself by translations}.

%%%%%%%%%%%%%%%%%%%%%%%%%%%%%%%%%%%%%%%%%%%%%%%%%%%%%%%%%%%%
\subsection{Flatness of Obata connection for exotic hypercomplex structures on a torus}
%%%%%%%%%%%%%%%%%%%%%%%%%%%%%%%%%%%%%%%%%%%%%%%%%%%%%%%%%%%%

%%%%%%%%%%%%%%%%%%%%%%%%%%%%%%%%%%%%%%%%%%%%%%%%%%%%%%%%%%%%
\theorem\label{_exotic_Obata_flat_Theorem_}
Let $(I,J,K)$ be a hypercomplex structure on a complex
torus $(M,I)$ and $\nabla$ the corresponding Obata connection.
Then $(M, I,J,K)$ is a quotient of $\H^n$ by an affine
action of $\Z^{4n}$.  In particular, $\nabla$ is flat.

\hfill

\pstep
Consider the affine structure on $(M,I)$ obtained from the
holomorphic structure as in \ref{_flat_on_torus_Remark_}.
 Since the fibers of the holomorphic projection
$\pi:\;\Tw(M) \arrow \C P^1$ are affine manifolds, 
the universal cover $\tilde{\pi}:  \Tw(\widetilde M) \to \C P^1$ is an affine
bundle.  Fixing a horizontal
section of $\pi$ is the same as fixing an origin for the
torus. This identifies $ \Tw(\widetilde M)$ with
$\Tot(\calo(1)^{2n})$ (the linearization of this affine
bundle is determined by \ref{_O(1)_normal_Remark_}). 

The anticomplex involution $\tau$ of $\Tw(M)$ clearly
respects the affine structure of the fibers and its lift
$\tilde{\tau}$ to ${\Tw}(\widetilde M)$ preserves the vector bundle structure.

\hfill

{\bf Step 2:} Consider the twistor data on 
$\Tot(\calo(1)^{2n})=\Tw(\widetilde M)$ which are compatible 
with the vector bundle structure, as in Step 1. Such twistor data are
uniquely determined by the vector bundle structure;
indeed, the space of sections of the twistor projection
has only one component, and the anticomplex involution
on a vector bundle is unique up to an automorphism
of $\calo(1)^{2n}$ which has constant coefficients,
hence is standard in appropriate coordinates. 
This implies that $\widetilde M$ is isomorphic to 
${\Bbb H}^n$ as a hypercomplex manifold.

\hfill

{\bf Step 3:}
The map $\Tw(\widetilde M) \to \Tw(M)$ is a morphism
of twistor data, locally an
isomorphism. Therefore,  the Obata connection of
$(M,I,J,K)$ is flat as well. Moreover, $\pi_1(\Tw(M))$
acts on ${\Tw}(\tilde M)$ preserving the fibers of
$\tilde{\pi}:\; \Tw(\widetilde M) \to\C P^1$. This action is affine on each fiber, 
because hypercomplex automorphisms of ${\Bbb H}^n$ 
are affine. Therefore,  $(M,I,J,K)$ is a quotient
of the flat affine hypercomplex manifold ${\Bbb H}^n$ by 
an affine action of $\Z^{4n}$.
\endproof

\hfill

\remark
If the holonomy of Obata connection
on $M$ is trivial (or just compact), it would
immediately follow that $M$ is a hyperk\"ahler torus.

\hfill

\corollary\label{_exotic_hc_torus_holonomy_Corollary_}
Let $(M,I,J,K)$ be a hypercomplex manifold, and $\nabla$ 
its Obata connection. Assume that $(M,I)$ is a compact
complex torus. Then $\nabla$ is flat and the hypercomplex
structure is exotic if and only if the holonomy
of $\nabla$ is non-trivial. \endproof

%%%%%%%%%%%%%%%%%%%%%%%%%%%%%%%%%%%%%%%%%%%%%%%%%%%%%%%%%%%%

\section{Flat affine structures on a torus}
\label{_flat_affine_on_torus_Section_}
%%%%%%%%%%%%%%%%%%%%%%%%%%%%%%%%%%%%%%%%%%%%%%%%%%%%%%%%%%%%

%%%%%%%%%%%%%%%%%%%%%%%%%%%%%%%%%%%%%%%%%%%%%%%%%%%%%%%%%%%%%%%%%%%%%%%%
\subsection{Non-standard flat affine structures on a torus}
%%%%%%%%%%%%%%%%%%%%%%%%%%%%%%%%%%%%%%%%%%%%%%%%%%%%%%%%%%%%%%%%%%%%%%%%

\remark
A flat affine structure on a torus $T$ is called {\bf
  standard} if its linear holonomy is trivial. It is
equivalent to say that $T$ is
isomorphic to $\R^n/ \Lambda$ as an affine manifold, where $\Lambda$ is an
integral lattice of rank $n$.

\hfill

\remark \label{_standard_flat_torus__via_holonomy_Remark_}
Let $(M, \nabla)$ be a standard flat affine torus. Then $\pi_1(M)$ acts on $\tilde M$
by translations, hence $\tilde M=\R^n$ and $M$ is
isomorphic to $\R^n/\Z^n$ with the standard flat affine
structure.

\hfill

\remark
In \cite{_Thurston_Sullivan:Charts_},
W. Thurston and D. Sullivan give examples of non-standard 
flat affine structures on a torus. In \cite{_baues:deformations_},
O. Baues discusses the deformation space of flat affine
structures on 2-torus.

\hfill

\example\label{_non_complete_flat_torus_Example_}
Consider the quotient
$M=\frac{\R^2\backslash 0}{\Z}$, where $\Z$ acts by homotheties.
The linear holonomy of $M$ is $\Z$ acting 
on $TM$ by homotheties. This is also an 
example of non complete affine manifold.

\hfill

\example\label{_exotic_flat_torus_main_Example_}
Consider the $\Z^2$-action $\rho$ on $\R^2$ generated by 
$(x, y) \to (x+1, y)$ and $(x, y) \to (x+y, y+1)$.
The projection to the second component maps
the quotient $X$ to $S^1$, with the fiber $S^1$,
hence $X$ is a compact manifold. It is diffeomorphic to a
torus $\R^2/\Z^2$ by construction. Its linear holonomy
is generated by $A(x,y):= (x+y, y)$.
Therefore this action defines a non-standard flat
affine structure on a torus.

%%%%%%%%%%%%%%%%%%%%%%%%%%%%%%%%%%%%%%%%%%%%%%%%%%%%%%%%%%%%%%%%%%%%%%%%%
\subsection{Fried-Goldman-Hirsch theorem}
%%%%%%%%%%%%%%%%%%%%%%%%%%%%%%%%%%%%%%%%%%%%%%%%%%%%%%%%%%%%%%%%%%%%%%%%%

The tools needed to show the non-existence of the exotic hypercomplex
structures on tori are provided by Fried, Goldman and
Hirsch, who laid out the groundwork in (\cite{_FGH:Affine_}).

From  \ref{_exotic_Obata_flat_Theorem_}, the following is clear:

\hfill

\proposition
Any hypercomplex torus $M$ is a complete affine manifold.

\proof Indeed, $M$ is a quotient of ${\Bbb H}^n$ by 
a discrete group preserving the affine structure. \endproof

\hfill

\conjecture \label{_Marcus_conjecture_}
(``Marcus conjecture'', \cite{_Abels:survey_})\\
A compact flat affine manifold is complete if and only if
it admits a parallel volume form.

\hfill

Fried, Goldman and Hirsch proved 
the Marcus conjecture, assuming that the affine holonomy
group is nilpotent.

\hfill

%%%%%%%%%%%%%%%%%%%%%%%%%%%%%%%%%%%%%%%%%%%%%%%%%%%%%%%%%%%%%
\theorem\label{_FGH_Theorem_}
Let $(M, \nabla)$ be a compact flat affine manifold 
with the affine holonomy group nilpotent.  Then the following
are equivalent:
\begin{description}
\item[(a)]  $(M, \nabla)$
is complete, 
 \item[(b)] $(M, \nabla)$ admits a paralell volume
form, and 
 \item[(c)]  its linear holonomy action is unipotent.
\end{description}

\proof
\cite[Theorem A]{_FGH:Affine_}.
\endproof

\hfill

\remark 
\ref{_FGH_Theorem_} can be applied to any
complete flat affine torus; it implies that its
linear holonomy is unipotent, unlike
\ref{_non_complete_flat_torus_Example_}.

\hfill

\theorem \label{unipotent}
Let $(M,I,J,K)$ a hypercomplex torus and
$\nabla$ its Obata connection. If $(M,I)$ admits a K\"{a}hler metric, then $(M,\nabla)$ satisfies (a)-(c) of Fried-Goldman-Hirsch theorem.

\hfill

\proof Since $(M,I)$ is K\"ahler, it is HKT (\cite{_Verbitsky:Kahler_HKT_}).
Since its canonical bundle is trivial and $(M,I,J,K)$ is
HKT, the Obata holonomy is contained in $SL(n, {\Bbb H})$
(this observation was essentially proven in
\cite{_Verbitsky:HKT_} and spelled out explicitly in
\cite{_Verbitsky:canoni_}). 
In other words, $\nabla$ fixes a volume form.
\endproof

%%%%%%%%%%%%%%%%%%%%%%%%%%%%%%%%%%%%%%%%%%%%%%%%%%%%%%%%%%%%%%%%%%%%%%%%
\subsection{Complete flat affine structures on real tori}
%%%%%%%%%%%%%%%%%%%%%%%%%%%%%%%%%%%%%%%%%%%%%%%%%%%%%%%%%%%%%%%%%%%%%%%%

We start from a classification of complete flat affine
structures on a torus.
Several arguments in the proofs rely on Maltsev's work,
originally contained in \cite{_Maltsev_}. We suggest to
consult \cite{_Raghunathan:discrete_subgroups_,_GO'H:nilpotent_} for a
comprehensive description.

\hfill

%%%%%%%%%%%%%%%%%%%%%%%%%%%%%%%%%%%%%%%%%%%%%%%%%%%%%%%%%%%%
\theorem\label{_real_tori_classified_Theorem_}
Let $(M,\nabla)$ be a flat affine structure on a
compact torus. Assume that $\nabla$ is complete.
Then 
\begin{description}
\item[(i)] 
its linear holonomy is unipotent. 
\item[(ii)]   
For some real basis $t_1, ..., t_{n}$ in $\R^n$, the 
action $\Gamma:=\pi_1(M)$ on $\R^n$,  considered 
as a universal cover of $M$, is generated by
$\tau_1,..., \tau_{n}$, with $\tau_i(x):= t_i +L_i(x)$.
Here $L_1, ..., L_{n}\in GL(n, \R)$ is a collection of commuting
unipotent matrices which satisfy
\begin{equation}\label{_commutators_affine_Equation_}
(L_i-\Id)(t_j)=(L_j-\Id)(t_i)
\end{equation}
 for any $i, j$.
\item[(iii)]
For any collection of commuting affine maps 
$\tau_1,..., \tau_{n}$ with $\tau_i(x):= t_i +L_i(x)$, 
with $t_1, ..., t_n$ a basis and 
all $L_i$ unipotent and satisfying \eqref{_commutators_affine_Equation_},
there exists a flat affine structure on a torus $M$ with $\tau_1,..., \tau_{n}$
generating the action of $\pi_1(M)$ on $\R^n$.
\end{description}

{\bf Proof of (i):}
	The linear holonomy of $\nabla$ is unipotent
by Fried-Goldman-Hirsch's \ref{_FGH_Theorem_},
because it is complete and $\pi_1(M)$ is nilpotent.

\hfill

\textbf{Proof of (ii), Step 1:}
Let us write the generators of $\Gamma=\pi_1(M)$
as $\tau_1,..., \tau_{n}$, with $\tau_i(x):= L_i(x)+t_i$,
where $L_1, ..., L_{n}$. Clearly, $L_i$ are commuting
and (by \ref{_FGH_Theorem_}) unipotent. 
The equation \eqref{_commutators_affine_Equation_}
follows because $\tau_i \tau_j (x) = L_i(t_j) + L_iL_j(x)+ t_i$,
which gives
\[
0=[\tau_i, \tau_j]=(L_i-\Id)(t_j)- (L_j-\Id)(t_i).
\]
To complete the proof of
(ii), it remains to show that 
all the $t_i$ are linearly independent.

\hfill

{\bf Step 2:} Let $G$ be the Maltsev completion of
$\Z^{n}=\pi_1(M)$, identified with the
Zariski closure of $\pi_1(M)$ in $\Aff(R^n)$ 
(\cite{_GO'H:nilpotent_}). 
Since $G$ is a nilpotent simply connected
Lie group, $\exp: \mathfrak{g} \to G$ is a
diffeomorphism. Moreover, it is possible to
find a basis $(X_i)$ of $\mathfrak{g}$ such that 
$\langle \exp(X_1), \dots, \exp(X_n)\rangle=\rho(\Z^n)\subset G$.
The canonical left invariant connection $\nabla$ on $G$ is both left and right invariant, since the left action is equal to
the right action. It is torsion free (since
$T(X,Y)=[X,Y]$) and flat, thus bestows $G$ with a flat
affine structure. In particular the inverse map $\log:\; G
\to \g$ coincides with the
development map $dev: G \to \R^{n}$, explicitly: $g\in G$ is
mapped to  $g\cdot 0\in \R^{n}$. Passing to the
differential 
$d_e(dev): \mathfrak{g} \to T_0\R^{n} $, a tangent vector $X\in
\mathfrak{g}$ is identified with $ \frac{d}{ds}\left(
\exp(sX)\cdot 0  \right)|_{s=0} \in \R^{n}  $. In particular,
setting $X_i=\log(\rho(z_i))$, we have
 \[
t_i=\exp(X_i)\cdot0=\frac{d}{ds}\left(   \exp(sX_i)\cdot 0  \right)\restrict{s=0} 
\]
i. e. under $d_e (dev) $ every $X_i$ is 
identified with $t_i$. Then they are all linearly independent.
This conlcudes the proof of (ii).

\hfill

{\bf Proof of (iii):}
Let $V=\R^n$.
We need to show that the action of the group $\langle \tau_1, ...,
\tau_n\rangle\subset \Aff(V)$ is free and cocompact. 
Let $V_0 \subset V_1 \subset ...\subset V$ be a filtration
on $V$ such that $L_i=\Id$ on its associated graded space $V_{gr}$
(such a filtration exists by Lie theorem, because all $L_i$ are unipotent
and commute). Since the constant terms $t_1,..., t_n$
are a basis in $V$, the group $\langle \tau_1,...,\tau_n\rangle$ acts on $V_{gr}$
as a lattice, and the quotient $V_{gr}/\Z^n$ is homeomorphic
to a torus. Then the action of
$\langle\tau_1,...,\tau_n\rangle$ on $V$ is free; indeed, 
if it has a fixed point, it would have a fixed point
on the associated graded space. Identifying $V$ with the
Maltsev completion of $\langle \tau_1,...,\tau_n\rangle$
as in Step 2 above, we obtain that $\langle
\tau_1,...,\tau_n\rangle$
is a lattice in a commutative Lie group, hence the
quotient $V/ \langle \tau_1,...,\tau_n\rangle$ is compact.
\endproof

\hfill

%%%%%%%%%%%%%%%%%%%%%%%%%%%%%%%%%%%%%%%%%%%%%%%%%%%%%%%%%%%%%%%%%%%%%%%%%
%\subsection{Complete flat affine structures on complex tori}
%%%%%%%%%%%%%%%%%%%%%%%%%%%%%%%%%%%%%%%%%%%%%%%%%%%%%%%%%%%%%%%%%%%%%%%%%

\remark\label{_3-tensor_symme_Remark_}
In assumptions of \ref{_real_tori_classified_Theorem_}, let $\Psi$ be the tensor
written in the basis $t_1, ..., t_{2n}$ as 
\[
\Psi(u,v):=\sum_{i,j,k=1}^{2n} (A_i)_j^k u^i v^j t_k,
\]
where $A_i=L_i-\Id$. Then $\Psi\in \Sym^{2} V^* \otimes V$.
If, in addition, the flat connection preserves a
complex or a hypercomplex structure on a torus,
this would mean that $A_i \in GL(n, \C)$
or  $A_i \in GL(n, {\Bbb H})$.

\hfill

\example    Let $W_1,W_2$ be complex
vector spaces. Take $V=W_1\oplus W_2$
and let $A_0\in \Hom_\C(W_1, W_2)$ be any 
complex linear map. Define $A\in \End(V)$ as $A\restrict{W_1}=A_0$, $A\restrict{W_2}=0$.
Write this matrix in a basis as $A:= \sum_{i,j}a_{j}^i t_i \otimes t^j$.
Take $\Psi:= \sum_{i,j}a_{j}^i t_i \otimes t_i \otimes t^j$.
This tensor is symmetric on first two variables
and complex linear on last two; also, $L:=-A+\Id$ is 
unipotent, defining a collection of 
complex linear affine transforms satisfying
the assumptions of \ref{_real_tori_classified_Theorem_} (iii).

\subsection{Exotic hypercomplex structures on tori do not  exist}
%%%%%%%%%%%%%%%%%%%%%%%%%%%%%%%%%%%%%%%%%%%%%%%%%%%%%%%%%%%%%%%%%%%%%%%%

\theorem\label{_no_hc_tori_exotic_Theorem_}
Exotic hypercomplex structures on tori do not exist.

\hfill

\proof
Let $(M,I,J,K,\nabla)$ be
a hypercomplex structure and Obata
connection on a complex torus.
By \ref{_exotic_Obata_flat_Theorem_},
$\nabla$ is flat and complete, hence 
\ref{_real_tori_classified_Theorem_} 
can be applied. Let
 $\Psi\in\Sym^2(V) \otimes V^*$ be the $(2,1)$-tensor obtained as 
in \ref{_3-tensor_symme_Remark_}, which now by hypothesis
commutes with $I,J,K$ on last two arguments.
This forces $\Psi$ to be zero, since
\[
K\Psi(x,y)=\Psi(IJx,y)=\Psi(Jx, Iy)=\Psi(x, JIy)=-K\Psi(x,y)
\]
and $K$ is a real automorphism of $V$.
\endproof

\hfill

\remark
The argument which proves
\ref{_no_hc_tori_exotic_Theorem_}
also proves the uniqueness of Obata connection.
Indeed, we prove that, given a quaternionic vector space, 
there are no  tensors $\Psi\in V^*\otimes V^*\otimes V$
which are symmetric in the first two variables and
${\Bbb H}$-linear on last two.
If $\nabla$ is a torsion free and quaternionic connection, then $\nabla=\nabla^{\text{Ob}} + \Psi$ and $\Psi \in \Sym^2{T^*M} \otimes TM$, commuting with $I,J,K$. Since $\Psi$ vanishes everywhere, $\nabla=\nabla^{\text{Ob}} $, as shown in \cite{_Obata_}.

\hfill

{\bf Achnowledgements:} 
We are grateful to Jorge Vitorio Pereira for the reference
to \cite{_Thurston_Sullivan:Charts_} and insightful comments, and to
Andrey Soldatenkov, Gueo Grantcharov and Anna Fino for their interest, 
suggestions and commentary.

\hfill

{\scriptsize

\noindent
{\sc Misha Verbitsky\\
 Instituto Nacional de Matem\'atica Pura e
	Aplicada (IMPA) \\ Estrada Dona Castorina, 110\\
	Jardim Bot\^anico, CEP 22460-320\\
	Rio de Janeiro, RJ - Brasil \\
{\tt verbit@impa.br}
}}

{\scriptsize
\noindent
{\sc Alberto Pipitone Federico\\
	Università degli Studi Tor Vergata
	 \\Via delle Ricerca Scientifica, 1 \\
	Roma, Italy \\
	{\tt pipitone@mat.uniroma2.it}
}}

\end{document}